\documentclass[a4paper,11pt,twoside]{article}
\usepackage{amsmath}
\usepackage{amsthm}
\usepackage{amsfonts, amssymb}
\usepackage{mathrsfs}
\usepackage[all]{xy}
\usepackage{latexsym}
\usepackage[dvips]{graphicx}

\title{\Large\sc{A Geometric Decomposition of Spaces into Cells of Different Types}}
\author{\sc{Gabriel Minian - Miguel Ottina}}

\date{}
\def\email#1{\vspace{1cm}\noindent E-mail address: {\sf #1} \bigskip}
\theoremstyle{plain}
\newtheorem{lemma}{Lemma}[section]
\newtheorem{prop}[lemma]{Proposition}
\newtheorem{theo}[lemma]{Theorem}
\newtheorem{coro}[lemma]{Corollary}
\newtheorem{theoint}[lemma]{Theorem}
\theoremstyle{remark}
\newtheorem{rem}[lemma]{Remark}

\theoremstyle{definition}
\newtheorem{definition}[lemma]{Definition}
\newtheorem{ex}[lemma]{Example}

\newtheorem{exs}[lemma]{Examples}

\newtheorem{imprem}[lemma]{Important remark}

\def\col{\mathrm{colim}\ }

\def\id{\mathrm{Id}}

\def\incl{\mathrm{inc}}

\def\cn{\mathrm{C}}

\begin{document}

\maketitle

\noindent{Departamento  de Matem\'atica\\
FCEyN, Universidad de Buenos Aires\\
Buenos Aires, Argentina.}



\begin{abstract}
\noindent We develop the theory of CW(A)-complexes, which generalizes the classical theory of CW-complexes, keeping
 the geometric intuition of J.H.C. Whitehead's original theory. We obtain this way generalizations of classical results, such
 as Whitehead Theorem, which allow a deeper insight in the homotopy properties of these spaces.
\end{abstract}



\section{Introduction}

\bigskip

It is well known that CW-complexes are spaces which are built up out of simple building blocks or \it cells. \rm 
In this case, 
balls are used as models for the cells and these are attached step by step using \it attaching maps\rm, which are defined in the 
boundary spheres of the balls. Since their introduction by J.H.C. Whitehead in the late fourties \cite{Whi},  CW-complexes 
have played an essential role in  geometry and topology. The combinatorial structure of these spaces allows 
the development of tools and results (e.g. simplicial and cellular aproximations, Whitehead Theorem, 
Homotopy excision, etc.) which lead to a deeper insight of their homotopy and homology properties.

\medskip

The main properties  of CW-complexes arise from the following two basic facts: (1) The $n$-ball $D^n$ is the 
topological (reduced) cone
of the $(n-1)$-sphere $S^{n-1}$ and (2) The $n$-sphere is the (reduced) $n$-suspension of the $0$-sphere $S^0$. 
For example, the homotopy 
extension properties of CW-complexes are deduced from (1), since the inclusion of the $(n-1)$-sphere in 
the $n$-disk is a closed cofibration. 
Item (2) is closely related to the definition of classical homotopy groups of spaces and it is used to prove results such as
Whitehead Theorem or Homotopy excision and in the construction of Eilenberg-MacLane spaces.
These two basic facts suggest also that one might replace the original \it core \rm $S^0$ by any other space $A$ and 
construct spaces built up out of cells of different \it shapes \rm or \it types \rm using suspensions and cones
of the base space $A$.

\medskip

The main purpose of this paper is to develop the theory of such spaces. More precisely, we define the notion of CW-complexes
of type $A$ (or CW(A)-spaces for short) generalizing CW-complexes (which constitute a special case of CW(A)-complexes, when
$A=S^0$). As in the classical case, we study these spaces from two different points of view: the constructive
 and the descriptive approachs. We use both points of view to prove  generalizations of classical results 
such as Whitehead Theorem and use these new results to study their homotopy properties. 

\medskip

Of course, some classical results are no longer true for general cores $A$. For 
example, the notion of dimension of a space (as a CW(A)-complex) is not always well defined. Recall that in the
classical case, the good definition of dimension is deduced from the famous Invariance of Dimension Theorem. By a 
similar argument,
we can prove that in particular cases (for example when the core $A$ is itself a finite dimensional CW-complex) 
the dimension of a CW(A)-complex is well defined. We study this and  other invariants and 
exhibit  many examples and counterexamples to clarify the main concepts.

\medskip

It is clear that, in general,  a topological space may admit many different decompositions into cells of different types. 
We study the relationship between such different decompositions. In particular, we obtain results such as the following.

\begin{theoint}
Let $A$ be a $CW(B)$-complex of finite dimension and let $X$ be a generalized $CW(A)$-complex. Then $X$ is a 
generalized $CW(B)$-complex. In particular, if $A$ is a standard finite dimensional CW-complex, then $X$ is a generalized CW-complex and
therefore it has the homotopy type of a CW-complex.
\end{theoint}

By a generalized complex we mean a space which is obtained by attaching cells in countable many steps, allowing cells of any dimension to be attached in any step.

\medskip

We also analyze the changing of the core $A$ by a core $B$ via a map $\alpha:A\to B$ and obtain the following result.

\begin{theoint}
Let $A$ and $B$ be pointed topological spaces with closed base points, let $X$ be a CW($A$) and 
let $\alpha:A\rightarrow B$ and $\beta:B\rightarrow A$ be continuous maps.
\begin{enumerate}
\item[i.] If $\beta\alpha=\id_A$, then there exists a CW($B$) $Y$  and maps $\varphi:X\rightarrow Y$ 
and $\psi:Y\rightarrow X$ such that $\psi\varphi=\id_X$.
\item[ii.] If $\beta$ is a homotopy equivalence, then there is a CW($B$) $Y$  and a homotopy 
equivalence $\varphi:X\rightarrow Y$.
\item[iii.] If $\beta\alpha=\id_A$ and $\alpha\beta\simeq\id_A$ then there exists a CW($B$) $Y$ 
 and maps $\varphi:X\rightarrow Y$ and $\psi:Y\rightarrow X$ such that $\psi\varphi=\id_X$ and $\varphi\psi=\id_X$.
\end{enumerate}
\end{theoint}

In particular, when the core $A$ is contractible, all CW($A$)-complexes are also contractible.

\medskip

Finally we start developing the homotopy theory of these spaces and obtain the following  generalization of Whitehead Theorem.

\begin{theoint}
Let $X$ and $Y$ be $CW(A)$-complexes and let $f:X\to Y$ be a continuous map. Then $f$ is a homotopy equivalence 
if and only if it is an $A$-weak equivalence.
\end{theoint}

We emphasize that our approach tries to keep the geometric intuition of Whitehead's original theory. There exist many 
generalizations of CW-complexes in the literature. We especially recommend Baues' generalization of complexes 
in Cofibration Categories \cite{Bau}. There is also a categorical approach to cell complexes by the first named author of this
paper \cite{Min}. The main advantage of the geometric point of view that we take in this article is that it allows the 
generalization of the most important classical results for CW-complexes and these new results can be applied in 
several concrete examples.


\bigskip

Throughout this paper, all spaces are assumed to be pointed spaces, all maps are pointed maps and  
homotopies are base-point preserving.

\section{The constructive approach and first results}

We denote by $\cn X$ the reduced cone of $X$ and by $\Sigma X$ its reduced suspension. Also, $S^n$ denotes the $n$-sphere and $D^n$ denotes the $n$-disk.

\bigskip

Let $A$ be a fixed pointed topological space.

\begin{definition}
We say that a (pointed) space \emph{$X$ is obtained from a (pointed) space $B$ by attaching an $n$-cell of type $A$} (or simply, an \emph{$A$-$n$-cell}) if there exists a pushout diagram
\begin{displaymath}
\xymatrix{\Sigma^{n-1}A \ar[r]^(.6){g} \ar[d]_i \ar@{}[rd]|{push} & B \ar[d] \\ \cn\Sigma^{n-1}A \ar[r]_(.6){f} & X}
\end{displaymath}
The \emph{$A$-cell} is the image of $f$. The map $g$ is the \emph{attaching map} of the cell, and $f$ is its \emph{characteristic map}.

We say that \emph{$X$ is obtained from $B$ by attaching a $0$-cell of type $A$} if $X=B\lor A$. 
\end{definition}

Note that attaching an $S^0$-$n$-cell is the same as attaching an $n$-cell in the usual sense, and that attaching an $S^m$-$n$-cell means attaching an $(m+n)$-cell in the usual sense. 

The reduced cone $\cn A$ of $A$ is obtained from $A$ by attaching an $A$-1-cell. In particular, $D^2$ is obtained from $D^1$ by attaching a $D^1$-1-cell. Also, the reduced suspension $\Sigma A$ can be obtained from the singleton $\ast$ by attaching an $A$-1-cell.

\bigskip

Of course, we can attach many $n$-cells at the same time by taking various copies of $\Sigma^{n-1}A$ and $\cn\Sigma^{n-1}A$.
\begin{displaymath}
\xymatrix{\underset{\alpha\in J}{\bigvee}\Sigma^{n-1}A \ar[r]^(.6){\underset{\alpha\in J}{+}g_\alpha} \ar[d]_i \ar@{}[rd]|{push} & B \ar[d] \\ \underset{\alpha\in J}{\bigvee}\cn\Sigma^{n-1}A \ar[r]_(.6){\underset{\alpha\in J}{+}f_\alpha} & X}
\end{displaymath}

\begin{definition}
A \emph{CW-structure with base $A$ on a space $X$}, or simply a \emph{CW($A$)-structure on $X$}, is a sequence of spaces $\ast=X^{-1}, X^0, X^1, \ldots,X^n,\ldots$ such that, for $n\in\mathbb{N}_0$, $X^n$ is obtained from $X^{n-1}$ by attaching $n$-cells of type $A$, and $X$ is the colimit of the diagram
$$\ast=X^{-1}\rightarrow X^0\rightarrow X^1 \rightarrow \ldots \rightarrow X^n \rightarrow \ldots$$
We call $X^n$ the \emph{$n$-skeleton} of $X$.

We say that the space $X$ is a \emph{CW($A$)-complex} (or simply a \emph{CW($A$)}), if it admits some CW($A$)-structure. In this case, the space $A$ will be called the \emph{core} or the \emph{base space} of the structure.
\end{definition}

\noindent
Note that a CW($A$) may admit many different structures of CW-complex with base $A$.

\medskip

\begin{exs} \ 
\begin{enumerate}
\item A CW($S^0$) is just a CW-complex and a CW($S^n$) is a CW-complex with no cells of dimension less than $n$.
\item The space $D^n$ admits several different CW($D^1$)-structures. For instance, we can take $X^r=D^{r+1}$ for $0 \leq r \leq n-1$ since $\cn D^r = D^{r+1}$. We may also take $X^0=\ldots = X^{n-2}=\ast$ and $X^{n-1}=D^n$ since there is a pushout
\begin{displaymath}
\xymatrix{\Sigma^{n-2}D^1=D^{n-1} \ar[r] \ar[d]_i \ar@{}[rd]|{push} & \ast \ar[d] \\ \cn\Sigma^{n-2}D^1=\cn D^{n-1} \ar[r] & \Sigma D^{n-1}=D^n}
\end{displaymath}
\end{enumerate}
\end{exs}

\bigskip

As in the classical case, instead of starting attaching cells from a base point $\ast$, we can start attaching cells on a pointed space $B$.

A \emph{relative CW($A$)-complex} is a pair $(X,B)$ such that $X$ is the colimit of a diagram
$$B=X^{-1}_B\rightarrow X^0_B\rightarrow X^1_B \rightarrow \ldots \rightarrow X^n_B \rightarrow \ldots$$
where $X^n_B$ is obtained from $X^{n-1}_B$ by attaching $n$-cells of type $A$.

\bigskip

It is clear that one can build a space $X$ by attaching cells (of some type $A$) without requiring them to be attached in such a way that their dimensions form an increasing sequence. That means, for example, that a 2-cell may be attached on a 5-cell.
In general, those spaces might not admit a CW($A$)-structure and they will be called \emph{generalized CW($A$)-complexes} (see \ref{defofgencwa}). If the core $A$ is itself a CW-complex, then a generalized CW($A$)-complex has the homotopy type of a CW-complex. This generalizes the well-known fact that a \emph{generalized CW-complex} has the homotopy type of a CW-complex.

Before we give the formal definition we show an example of a generalized CW-complex which is not a CW-complex.

\begin{ex} \label{remcwavsgen}
We build $X$ as follows. We start with a $0$-cell and we attach a $1$-cell by the identity map obtaining the interval $[-1;1]$. We regard $1$ as the base point. Now, for each $n\in\mathbb{N}$ we define $g_n:S^0\rightarrow [-1,1]$ by $g_n(1)=1$, $g_n(-1)=1/n$. We attach 1-cells by the maps $g_n$. This space $X$ is an example of a generalized CW-complex (with core $S^0$).

It is not hard to verify that it is not a CW-complex. To prove this, note that the points of the form $1/n$ must be 0-cells by a dimension argument, but they also have a cluster point at 0.
\end{ex}

\begin{definition} \label{defofgencwa}
We say that \emph{$X$ is obtained from $B$ by attaching cells} (of different dimensions) \emph{of type $A$} if there is a pushout
\begin{displaymath}
\xymatrix{\underset{\alpha\in J}{\bigvee}\Sigma^{n_\alpha-1}A \ar[r]^(.6){\underset{\alpha\in J}{+}g_\alpha} \ar[d]_i \ar@{}[rd]|{push} & B \ar[d] \\ (\underset{\alpha\in J_0}{\bigvee}A)\lor(\underset{\alpha\in J}{\bigvee}\cn\Sigma^{n_\alpha-1}A) \ar[r]_(.75){\underset{\alpha\in J}{+}f_\alpha} & X}
\end{displaymath}
where $n_\alpha\in\mathbb{N}$ for all $\alpha\in J$.
We say that $X$ is a \emph{generalized CW($A$)-complex} if $X$ is the colimit of a diagram
$$\ast=X^0 \rightarrow X^1 \rightarrow X^2 \rightarrow \ldots \rightarrow X^n \rightarrow \ldots$$
where $X^n$ is obtained from $X^{n-1}$ by attaching cells (of different dimensions) of type $A$.

We call $X^n$ the \emph{$n$-th layer} of $X$.

One can also define generalized relative CW($A$)-complexes in the obvious way.
\end{definition}

For standard CW-complexes, by the classical \emph{Invariance of Dimension Theorem}, one can prove that the notion of dimension is well defined. Any two different structures of a CW-complex must have the same dimension.

For a general core $A$ this is no longer true. However, we shall prove later that for particular cases (for example when $A$ is a finite dimensional CW-complex) the notion of dimension of a CW($A$)-complex is well defined.

\begin{definition}
Let $X$ be a CW($A$). We consider $X$ endowed with a particular CW($A$)-structure $\mathcal{K}$.
We say that the \emph{dimension} of $\mathcal{K}$ is $n$ if $X^n=X$ and $X^{n-1}\neq X$, and we write $\dim(\mathcal{K})=n$. We say that $\mathcal{K}$ is \emph{finite dimensional} if $\dim(\mathcal{K})=n$ for some $n\in\mathbb{N}_0$. 
\end{definition}

\begin{imprem}
A CW($A$) may admit different CW($A$)-structures with different dimensions. For example, let $A=\underset{n\in\mathbb{N}}{\bigvee}S^n$ and let $X=\underset{j\in\mathbb{N}}{\bigvee}A$. Then $X$ has a zero-dimensional CW($A$) structure. But we can see $X=(\underset{j\in\mathbb{N}}{\bigvee}A)\lor\Sigma A$, which induces a 1-dimensional structure. Note that $\underset{j\in\mathbb{N}}{\bigvee}A=(\underset{j\in\mathbb{N}}{\bigvee}A)\lor\Sigma A$ since both spaces consist of countably many copies of $S^n$ for each $n\in\mathbb{N}$.

\bigskip

Another example is the following. It is easy to see that if $B$ is a topological space with the indiscrete topology then its reduced cone and suspension also have the indiscrete topology. So, let $A$ be an indiscrete topological space with $1\leq \#A\leq c$. If $A$ is just a point then its reduced cone and suspension are also singletons, so $\ast$ can be given a CW($\ast$) structure of any dimension. If $\#A\geq 2$ then $\#(\Sigma^n A)=c$ for all $n$, and $\Sigma^n A$ are all indiscrete spaces. Since they have all the same cardinality and they are indiscrete then all of them are homeomorphic. But each $\Sigma^n A$ has an obvious CW($A$) structure of dimension $n$. Thus, the homeomorphisms between $\Sigma^n A$ and $\Sigma^m A$, for all $m$, allow us to give $\Sigma^n A$ a CW($A$) structure of any dimension (greater than zero).
\end{imprem}

Given a CW($A$)-complex $X$, we define the \emph{boundary} of an $n$-cell $e^n$ by $\overset{\bullet}{e^n}=e^n\cap X^{n-1}$ and the \emph{interior} of $e^n$ by $\overset{\circ}{e^n}=e^n-\overset{\bullet}{e^n}$. 

A cell $e^m_\beta$ is called an \emph{immediate face of $e^n_\alpha$} if $\overset{\circ}{e^m_\beta}\cap e^n_\alpha\neq\varnothing$, and a cell $e^m_\beta$ is called a \emph{face of $e^n_\alpha$} if there exists a finite sequence of cells $$e^m_\beta=e^{m_0}_{\beta_0}, e^{m_1}_{\beta_1}, e^{m_2}_{\beta_2}, \ldots, e^{m_k}_{\beta_k}=e^n_\alpha$$
such that $e^{m_j}_{\beta_j}$ is an immediate face of $e^{m_{j+1}}_{\beta_{j+1}}$ for $0\leq j<k$.

Finally, we call a cell \emph{principal} if it is not a face of any other cell.

\begin{rem}
Note that $\overset{\circ}{e^n_\alpha}\cap\overset{\circ}{e^m_\beta}\neq\varnothing$ if and only if $n=m$, $\alpha=\beta$. Thus, if $e^m_\beta$ is a face of $e^n_\alpha$ and $e^m_\beta \neq e^n_\alpha$ then $m<n$.
\end{rem}

\bigskip

As in the classical case, we can define subcomplexes and cellular maps in the obvious way.

\begin{rem}
If $X$ is a CW($A$), then $X=\underset{n,\alpha}{\bigcup}\overset{\circ}{e^n_\alpha}$.
\end{rem}


\begin{prop}
Let $X$ be a CW($A$) and suppose that the base point of $A$ is closed in $A$. Then the interiors of the $n$-cells are open in the $n$-skeleton. In particular, $X^{n-1}$ is a closed subspace of $X^n$.
\end{prop}

\begin{proof}
For $n=-1$ and $n=0$ it is clear. Let $n\geq 1$. We have a pushout diagram
\begin{displaymath}
\xymatrix{\underset{\alpha\in J}{\bigvee}\Sigma^{n-1}A \ar[r]^(.6){\underset{\alpha\in J}{+}g_\alpha} \ar[d]_i \ar@{}[rd]|{push} & X^{n-1} \ar[d] \\ \underset{\alpha\in J}{\bigvee}\cn\Sigma^{n-1}A \ar[r]_(.4){\underset{\alpha\in J}{+}f_\alpha} & X^n=X^{n-1}\cup\underset{\alpha}{\bigcup}e^n_\alpha}
\end{displaymath}
Consider a cell $e^n_\beta$. In order to verify that $\overset{\circ}{e^n_\beta}$ is open in $X^n$ we have to prove that $(+f_\beta)^{-1}(\overset{\circ}{e^n_\beta})$ is open in $\underset{\alpha\in J}{\bigvee}\cn\Sigma^{n-1}A$. Since $(+f_\beta)^{-1}(\overset{\circ}{e^n_\beta})=\cn\Sigma^{n-1}A-\Sigma^{n-1}A$ is open in $\cn\Sigma^{n-1}A$, then $\overset{\circ}{e^n_\beta}$ is open in $X^n$.
\end{proof}

\begin{prop} \label{gooddefofdim1}
Let $A$ be a finite dimensional CW-complex, $A\neq\ast$, and let $X$ be a CW($A$). Let $\mathcal{K}$ and $\mathcal{K}'$ be CW($A$)-structures in $X$ and let $n,m\in\mathbb{N}_0\cup\{\infty\}$ denote their dimensions. Then $n=m$.
\end{prop}

\begin{proof}
We suppose first that $\mathcal{K}$ and $\mathcal{K}'$ are finite dimensional and $n\geq m$. 

Let $k=\dim(A)$ and let $e^n_\alpha$ be an $n$-cell of $\mathcal{K}$. We have a homeomorphism $\overset{\circ}{e^n_\alpha}\simeq\cn\Sigma^{n-1}A-\Sigma^{n-1}A$, and $\overset{\circ}{e^n_\alpha}$ is open in $X$. Let $e$ be a cell of maximum dimension of the CW-complex $\cn\Sigma^{n-1}A$ and let $U=\overset{\circ}{e}$. Thus $U$ is open in $X$ and homeomorphic to $\overset{\circ}{D^{n+k}}$.

Now, $U$ intersects some interiors of cells of type $A$ of $\mathcal{K}'$. Let $e_0$ be one of those cells with maximum dimension. Suppose $e_0$ is an $m'$-cell, with $m'\leq m$. Then $\overset{\circ}{e_0}$ is open in the $m'$-skeleton of $X$ with the $\mathcal{K}'$ structure. It is not hard to see that $V=U\cap \overset{\circ}{e_0}$ is open in $U$, extending $\overset{\circ}{e_0}$ to an open subset of $X$ as in \ref{compintersfinitosint} below.

In a similar way, $\overset{\circ}{e_0}\simeq\cn\Sigma^{m'-1}A-\Sigma^{m'-1}A$, and $V$ meets some interiors of cells of the CW-complex $\cn\Sigma^{m'-1}A$. We take $e_1$ a cell (of type $S^0$) of maximum dimension among those cells and we denote $k'=\dim(e_1)$. Then $\overset{\circ}{e_1}$ is homeomorphic to $\overset{\circ}{D^{k'}}$. 
Let $W=V\cap\overset{\circ}{e_1}$. One can check that $W$ is open in $\overset{\circ}{e_1}\simeq\overset{\circ}{D^{k'}}$ and that it is also open in $U\simeq\overset{\circ}{D^{n+k}}$.

By the invariance of dimension theorem, $n+k=k'$, but also $k'\leq m+k \leq n+k$. Thus $n=m$.

\medskip

It remains to be shown that if $m=\infty$ then $n=\infty$. Suppose that $m=\infty$ and $n\neq\infty$. Let $k=\dim(A)$. We choose $e^l$ an $l$-cell of $\mathcal{K}'$ with $l>n+k$. Then $\overset{\circ}{e^l}$ is open in the $l$-skeleton $(\mathcal{K}')^l$. As in the proof of \ref{compintersfinitosint} below, we can extend $\overset{\circ}{e^l}$ to an open subset $U$ of $X$ with $U\cap (\mathcal{K}')^{l-1}=\varnothing$. Now we take a cell $e_1$ of $\mathcal{K}$ such that $\overset{\circ}{e_1}\cap U\neq \varnothing$ and with the property of being of maximum dimension among the cells of $\mathcal{K}$ whose interior meets $U$. Let $r=\dim(e_1)$. We have that $U\subseteq \mathcal{K}^r$. As before, we extend $\overset{\circ}{e_1}$ to an open subset $V$ of $X$ with $V\cap \mathcal{K}^{r-1}=\varnothing$, $V\cap \mathcal{K}^r=\overset{\circ}{e_1}$. So $U\cap\overset{\circ}{e_1}=U\cap V$ is open in $X$. Proceeding analogously, since $\overset{\circ}{e_1}\simeq\cn\Sigma^{r-1}A-\Sigma^{r-1}A$, we can choose a cell $e_2$ of $e_1$ (of type $S^0$) with maximum dimension such that $W=\overset{\circ}{e_2}\cap(U\cap\overset{\circ}{e_1})\neq\varnothing$. Again, $W$ is open in $X$. Let $s=\dim{e_2}$. So $W$ is open in $\overset{\circ}{e_2}\simeq \overset{\circ}{D^s}$ and $s\leq r+k\leq n+k <l$. On the other hand, $W$ must meet the interior of some cell of type $S^0$ belonging to one of the cells of $\mathcal{K}'$ with dimension greater than or equal to $l$ (since $U\cap (\mathcal{K}')^{l-1}=\varnothing$). So, a subset of $W$ is homeomorphic to an open set of $\overset{\circ}{D^q}$ with $q\geq l$, a contradiction.
\end{proof}

Recall that a topological space $Y$ is T1 if the points are closed in $X$.

\begin{prop} \label{compintersfinitosint}
Let $A$ be a pointed T1 topological space, let $X$ be a CW($A$) and $K\subseteq X$ a compact subspace. Then $K$ meets only a finite number of interiors of cells.
\end{prop}

\begin{proof}
Let $\Lambda=\{\alpha /\ K\cap \overset{\circ}{e^{n_\alpha}_\alpha}\neq\varnothing\}$. For each $\alpha\in\Lambda$ choose $x_\alpha\in K\cap\overset{\circ}{e^{n_\alpha}_\alpha}$. We want to show that for any $\alpha\in\Lambda$ there exists an open subspace $U_\alpha\subseteq X$ such that $U_\alpha\supseteq\overset{\circ}{e^{n_\alpha}_\alpha}$ and $x_\beta\notin U_\alpha$ for any $\beta\neq\alpha$. 

For each $n$, let $J_n$ be the index set of the $n$-cells. We denote by $g_\alpha^n$ the attaching map of $e_\alpha^n$ and by $f_\alpha^n$ its characteristic map.

Fix $\beta\in\Lambda$. Take $U_1=\overset{\circ}{e^{n_\beta}_\beta}$, which is open in $X^{n_\beta}$. If $n_\beta=-1$, we take $U_2=(\underset{\alpha\in J_0\cap\Lambda}{\bigvee}A-\{x_\alpha\})\lor(\underset{\alpha\in J_0-\Lambda}{\bigvee}A)$, which is open in the $0$-skeleton.

Now, for $n_\beta+n-1\geq 1$ we construct inductively open subspaces $U_n$ of $X^{n_\beta+n-1}$ with $U_{n-1}\subseteq U_n$, $U_n \cap X^{n_\beta+n-2}= U_{n-1}$ and such that $x_\alpha\notin U_n$ if $\alpha\neq\beta$.

If the base point $a_0\notin U_{n-1}$, we take $$U_n=U_{n-1}\cup\underset{\alpha\in J_{n_\beta+n-1} }{\bigcup}f_\alpha^{n_\alpha}( (g_\alpha^{n_\alpha})^{-1}(U_{n-1})\times (1-\varepsilon_\alpha,1])$$ with $0<\varepsilon_\alpha<1$ chosen in such a way that $x_\alpha\notin U_n$ if $\alpha\neq\beta$. Note that $U_n$ is open in $X^{n_\beta+n-1}$.

If $a_0\in U_{n-1}$ we take $$U_n=U_{n-1}\cup\underset{\alpha\in J_{n_\beta+n-1} }{\bigcup}f_\alpha^{n_\alpha}(( (g_\alpha^{n_\alpha})^{-1}(U_{n-1})\times (1-\varepsilon_\alpha,1])\cup (W_{x_\alpha}\times I) \cup (\Sigma^{n_\beta+n-1}A\times [0,\varepsilon'_\alpha)) )$$
with $W_{x_\alpha}= V_{x_\alpha}\cap(g_\alpha^{n_\alpha})^{-1}(U_{n-1})$, where $V_{x_\alpha}\subseteq\Sigma^{n_\beta+n-1}A$ is an open neighbourhood of the base point not containing $x_\alpha'$ (where $x_\alpha=f_\alpha^{n_\alpha}(x_\alpha',t_\alpha)$), and $0<\varepsilon_\alpha<1$, $0<\varepsilon'_\alpha<1$, chosen in such a way that $x_\alpha\notin U_n$ if $\alpha\neq\beta$. Note that $U_n$ is open in $X^{n_\beta+n-1}$.

We set $U_\beta=\underset{n\in\mathbb{N}}{\bigcup}U_n$.
Thus $K\subseteq\underset{\alpha\in\Lambda}{\bigcup} \overset{\circ}{e^{n_\alpha}_\alpha} \subseteq \underset{\alpha\in\Lambda}{\bigcup}U_\alpha$, and $x_\alpha\notin U_\beta$ if $\alpha\neq\beta$. Since $\{U_\alpha\}_{\alpha\in\Lambda}$ is an open covering of $K$ which does not admit a proper subcovering, $\Lambda$ must be finite.
\end{proof}

\begin{lemma}
Let $A$ and $B$ be Hausdorff spaces and suppose $X$ is obtained from $B$ by attaching cells of type $A$. Then $X$ is Hausdorff.
\end{lemma}

\begin{proof}
Let $x,y\in X$. If $x,y$ lie in the interior of some cell, then it is easy to choose the open neighbourhoods. If one of them belongs to $B$ and the other to the interior of a cell, let's say $x\in e^{n_\alpha}_\alpha$, we work as in the previous proof. Explicitly, if $x=f_\alpha(a,t)$ with $a\in\Sigma^{n_\alpha-1}A$, $t\in I$ then we take $U'\subseteq \Sigma^{n_\alpha-1}A$ open set such that $a\in U'$ and $a_0\notin U'$, where $a_0$ is the basepoint of $\Sigma^{n_\alpha-1}A$. We define $U=f_\alpha(U'\times (t/2,(1+t)/2))$, and $V=X-f_\alpha(\overline{U'}\times [t/2,(1+t)/2])$.

If $x,y\in B$, since $B$ is Hausdorff there exist $U',V'\subseteq B$ open disjoint sets such that $x\in U'$ and $y\in V'$. However, $U'$ and $V'$ need not be open in $X$. Suppose first that $x$, $y$ are both different from the base point. So we may suppose that neither $U'$ nor $V'$ contain the base point. We take $$U=U'\cup\underset{\alpha\in J}{\bigcup}f_\alpha( (g_\alpha)^{-1}(U')\times (1/2;1])$$
$$V=V'\cup\underset{\alpha\in J}{\bigcup}f_\alpha( (g_\alpha)^{-1}(V')\times (1/2;1])$$
If $x$ is the base point then we take 
$$U=U'\cup\underset{\alpha\in J}{\bigcup}f_\alpha( ((g_\alpha)^{-1}(U')\times I)\cup(\Sigma^{n_\alpha-1}A\times [0;1/2)))$$

%
\end{proof}

\begin{prop}
Let $A$ be a Hausdorff space and let $X$ be a CW($A$). Then $X$ is a Hausdorff space.
\end{prop}

\begin{proof}
By the previous lemma and induction we have that $X^n$ is a Hausdorff space for all $n\geq -1$. Given $x,y\in X$, choose $m\in\mathbb{N}$ such that $x,y\in X^m$. As $X^m$ is a Hausdorff space, there exist disjoint sets $U_0$ and $V_0$, which are open in $X^m$, such that $x\in U_0$ and $y\in V_0$. Proceeding in a similar way as we did in the previous results we construct inductively sets $U_k$, $V_k$ for $k\in\mathbb{N}$ such that $U_k,V_k\subseteq X^{m+k}$ are open sets, $U_k\cap V_k=\varnothing$, $U_k\cap X^{m+k-1}=U_{k-1}$ and $V_k\cap X^{m+k-1} =V_{k-1}$ for all $k\in\mathbb{N}$. We take $U=\bigcup U_k$, $V=\bigcup V_k$.
\end{proof}

\begin{rem} \label{remweaktopo}
Let $X$ be a CW($A$) and $S\subseteq X$ a subspace. Then $S$ is closed in $X$ if and only if $S\cap e^n_\alpha$ is closed in $e^n_\alpha$ for all $n$, $\alpha$.
\end{rem}

\begin{lemma} \label{cwpaste}
Let $X$, $Y$ be CW($A$)'s, $B\subseteq X$ a subcomplex, and $f:B\rightarrow Y$ a cellular map. Then the pushout
\begin{displaymath}
\xymatrix{B \ar[r]^f \ar[d]_i & Y \ar[d] \\ X \ar[r] \ar@{}[ru]|{push} & X\underset{B}{\cup}Y}
\end{displaymath}
is a CW($A$).
\end{lemma}

\begin{proof}

We denote by $\{e^n_{X,\alpha}\}_{\alpha\in J_n}$ the $n$-cells (of type $A$) of the relative CW($A$)-complex $(X,B)$ and by $\{e^n_{Y,\alpha}\}_{\alpha\in J'_n}$ the $n$-cells of $Y$. We will construct $X\underset{B}{\cup}Y$ attaching the cells of $Y$ with the same attaching maps and at the same time we will attach the cells of $(X,B)$ using the map $f:B\rightarrow Y$.

Let $J_0''=J_0\cup J_0'$ and $Z^0=\underset{\alpha\in J_0''}{\bigvee}A$. We define $f_0:X^0\rightarrow Z^0$ by $f_0|_{B^0}=f|_{B^0}$ and $f_0|_{\underset{\alpha\in J_0}{\bigcup} e^0_{X,\alpha}}$ the inclusion.

Suppose that $Z^{n-1}$ and $f_{n-1}:X^{n-1}\rightarrow Z^{n-1}$ with $f_{n-1}|_{B^{n-1}}=f$ are defined. We define $Z^n$ by the following pushout.
\begin{displaymath}
\xymatrix{\underset{\alpha\in J''_n}{\bigvee}\Sigma^{n-1}A \ar[r]^(.6){\underset{\alpha\in J''_n}{+}g''_\alpha} \ar[d]_i \ar@{}[rd]|{push} & Z^{n-1} \ar[d]^{i_{n-1}} \\ \underset{\alpha\in J''_n}{\bigvee}\cn\Sigma^{n-1}A \ar[r]_(.6){\underset{\alpha\in J''_n}{+}f''_\alpha} & Z^n}
\end{displaymath}
where $J''_n=J_n\cup J'_n$ and
\begin{displaymath}
g''_\alpha=\left\{\begin{array}{ll} f_{n-1}\circ g_\alpha & \textrm{if } \alpha \in J_n \\ g'_\alpha & \textrm{if } \alpha \in J'_n \end{array} \right.
\end{displaymath}
where $g_\alpha$ and $g_\alpha'$ are the attaching maps.
We define $f_n:X^n\rightarrow Z^n$ by $f_n|_{B^n}=f|_{B^n}$, $f_n|_{X^{n-1}}=f_{n-1}$ and $f_n|_{\underset{\alpha\in J_n}{\bigcup} e^n_{X,\alpha}}=f''_\alpha$ (i.e. $f_n(f_\alpha(x))=f''_\alpha(x)$). Note that $f_n$ is well defined.

Let $Z$ be the colimit of the $Z^n$. By construction it is not difficult to verify that $Z$ satisfies the universal property of the pushout.
\end{proof}

\begin{coro}
Let $X$ be a CW($A$) and $B\subseteq X$ a subcomplex. Then $X/B$ is a CW($A$).
\end{coro}

\begin{theo} \label{conoysuspcwa}
Let $X$ be a CW($A$). Then the reduced cone $\cn X$ and the reduced suspension $\Sigma X$ are CW($A$)'s. Moreover, $X$ is a subcomplex of both of them.
\end{theo}

\begin{proof}
By the previous lemma, it suffices to prove the result for $\cn X$.

Let $e^n_\alpha$ be the $n$-cells of $X$ and, for each $n$, let $J_n$ be the index set of the $n$-cells. We denote by $g^n_\alpha$ the attaching maps and  by $f^n_\alpha$ the characteristic maps. Let $i_{n-1}:X^{n-1}\rightarrow X^n$ be the inclusions. We construct $Y=\cn X$ as follows.


Let $Y^0=\underset{\alpha\in J_0}{\bigvee}A=X^0$.

We construct $Y^1$ from $Y^0$ and from the $0$-cells and the $1$-cells of $X$ by the pushout
\begin{displaymath}
\xymatrix{\underset{\alpha\in J'_1}{\bigvee}A \ar[r]^(.6){\underset{\alpha\in J'_1}{+}g'_\alpha} \ar[d]_i \ar@{}[rd]|{push} & Y^0 \ar[d]^{i_0'} \\ \underset{\alpha\in J'_1}{\bigvee}\cn A \ar[r]_(.6){\underset{\alpha\in J'_1}{+}f'_\alpha} & Y^1}
\end{displaymath}

where $J'_1=J_0 \sqcup J_1$. The maps $g'_\alpha$, for $\alpha\in J_1'$, are defined as
\begin{displaymath}
g'_\alpha=\left\{\begin{array}{ll} i_\alpha & \textrm{if } \alpha \in J_0 \\ g_\alpha & \textrm{if } \alpha \in J_1\end{array} \right.
\end{displaymath}
and $i_\alpha:A\rightarrow \underset{\alpha\in J_0}{\bigvee}A$ is the inclusion of $A$ in the $\alpha$-th copy. Note that $X^1$ is a subcomplex of $Y^1$.

Note also that the 1-cells of $Y$ are divided into two sets. The ones with $\alpha\in J_1$ are the 1-cells of $X$, and the others are the cone of the 0-cells of $X$.

Inductively, suppose we have constructed $Y^{n-1}$. We define $Y^n$ as the pushout
\begin{displaymath}
\xymatrix{\underset{\alpha\in J'_n}{\bigvee}\Sigma^{n-1}A \ar[r]^(.6){\underset{\alpha\in J'_n}{+}g'_\alpha} \ar[d]_i \ar@{}[rd]|{push} & Y^{n-1} \ar[d]^{i_{n-1}'} \\ \underset{\alpha\in J'_n}{\bigvee}\cn\Sigma^{n-1}A \ar[r]_(.6){\underset{\alpha\in J'_n}{+}f'_\alpha} & Y^n}
\end{displaymath}
where $J'_n=J_{n-1} \sqcup J_n$ and
\begin{displaymath}
g'_\alpha=\left\{\begin{array}{ll} g_\alpha & \textrm{for } \alpha \in J_n \\ f_\alpha\cup \cn g_\alpha & \textrm{for } \alpha \in J_{n-1} \ .\end{array} \right.
\end{displaymath}

We prove now that $Y^n=\cn X^{n-1}\cup\underset{\alpha}{\bigcup}e^n_\alpha$. We have the following commutative diagram.

\begin{displaymath} \label{eq:diagram}
\xymatrix@C=40pt{\underset{\alpha\in J'_n}{\bigvee}\Sigma^{n-1}A \ar[r]^(.4){(\underset{\alpha\in J_{n-1}}{+}g'_\alpha) \lor \id} \ar[d]_{\underset{\alpha\in J'_n}{\bigvee}i} & Y^{n-1}\lor\underset{\alpha\in J_n}{\bigvee}\Sigma^{n-1}A \ar[d]^{i_{n-1}\lor\underset{\alpha\in J_n}{\bigvee}i} \ar[r]^(.55){\id+(\underset{\alpha\in J_n}{+}g'_\alpha )} \ar@{}[rd]|{push} & Y^{n-1} \ar[d] \\ \underset{\alpha\in J'_n}{\bigvee}\cn\Sigma^{n-1}A \ar[r]_(.4){(\underset{\alpha\in J_{n-1}}{+}f'_\alpha) \lor \id} & \cn X^{n-1}\lor\underset{\alpha\in J_n}{\bigvee}\cn\Sigma^{n-1}A \ar[r]_(.55){\id+(\underset{\alpha\in J_n}{+}f'_\alpha)} & \cn X^{n-1} \cup\underset{\alpha}{\bigcup}e^n_\alpha}
\end{displaymath}

The right square is clearly a pushout. To prove that the left square is also a pushout it suffices to verify that the following is also a pushout.

\begin{displaymath}
\xymatrix@C=40pt{\underset{\alpha\in J_{n-1}}{\bigvee}\Sigma^{n-1}A \ar[r]^(.35){\underset{\alpha\in J_{n-1}}{+}g'_\alpha} \ar[d]_{\underset{\alpha\in J_{n-1}}{\bigvee}i} & Y^{n-1}=\cn X^{n-2}\cup \underset{\alpha\in J_{n-1}}{\bigcup}e^{n-1}_\alpha \ar[d]^{\textrm{inc}} \\ \underset{\alpha\in J_{n-1}}{\bigvee}\cn\Sigma^{n-1}A \ar[r]_(.6){\underset{\alpha\in J_{n-1}}{+}f'_\alpha} & \cn X^{n-1}}
\end{displaymath}

For simplicity, we will prove this in the case that there is only one $A$-(n-1)-cell. Let
\begin{displaymath}
\begin{array}{l}
j:\Sigma^{n-1}A\rightarrow \cn\Sigma^{n-1}A \\
i_1:\cn(\Sigma^{n-1}A)\times\{1\}\rightarrow\cn\cn\Sigma^{n-1}A \\ i_2:(\Sigma^{n-1}A)\times\{1\}\times I/\sim \rightarrow\cn\cn\Sigma^{n-1}A \\ i:\Sigma^{n}A=\cn\Sigma^{n-1}A\underset{A}{\cup}\cn\Sigma^{n-1}A\rightarrow\cn\Sigma^n A
\end{array}
\end{displaymath}
be the corresponding inclusions.

Let $\varphi:\cn\cn(\Sigma^{n-1}A)\rightarrow\cn\Sigma(\Sigma^{n-1}A)$ be a homeomorphism, such that $\varphi^{-1}i=i_1+i_2$. Note that $\cn j=i_2$.
There are pushout diagrams
\begin{displaymath}
\xymatrix{\Sigma^{n-1}A \ar[r]^{g} \ar[d]_{j} \ar@{}[rd]|{push} & X^{n-1} \ar[d]^{\textrm{inc}} \\ \cn\Sigma^{n-1}A \ar[r]_(.4){f} & X^{n}=X^{n-1}\cup e^{n}}
\qquad
\xymatrix{\cn\Sigma^{n-1}A \ar[r]^{\cn g} \ar[d]_{\cn j =i_2} \ar@{}[rd]|{push} & \cn X^{n-1} \ar[d]^{\cn\incl} \\ \cn\cn\Sigma^{n-1}A \ar[r]_(.55){\cn f} & \cn X^{n}}
\end{displaymath}

It is not hard to check that the diagram

\begin{displaymath}
\xymatrix@C=40pt{\Sigma^{n}A=\cn\Sigma^{n-1}A\underset{A}{\cup}\cn\Sigma^{n-1}A \ar[r]^(.6){f+\cn g} \ar[d]_{i} & \cn X^{n-1}\cup e^{n} \ar[d]^{\textrm{inc}} \\ \cn\Sigma^{n}A \ar[r]_{(\cn f) \varphi^{-1}} & \cn X^{n}}
\end{displaymath}

satisfies the universal property of pushouts.

%
%

Now we take $Y$ to be the colimit of $Y^n$, which satisfies the desired properties.
\end{proof}

\begin{rem} \ 
\begin{enumerate}

\item The standard proof of the previous theorem for a CW-complex $X$ uses the fact that $X\times I$ is also a CW-complex. For general cores $A$, it is not always true that $X\times I$ is a CW($A$)-complex when $X$ is.
%

\item It is easy to see that if $X$ is a CW($A$), then $\Sigma X$ is a CW($A$). Just apply the $\Sigma$ functor to each of the pushout diagrams used to construct $X$. In this way we give $\Sigma X$ a CW($A$) structure in which each of the cells is the reduced suspension of a cell of $X$. This is a simple and interesting structure. However, it does not have the property of having $X$ as a subcomplex.
\end{enumerate}
\end{rem}

\begin{lemma} \label{cwrelpaste}
Let $A$ be a topological space and let $(X,B)$ be a relative CW($A$) (resp. a generalized relative CW($A$)). Let $Y$ be a topological space, and let $f:B\rightarrow Y$ be a continuous map. We consider the pushout diagram
\begin{displaymath}
\xymatrix{B \ar[r]^f \ar[d]_i & Y \ar[d] \\ X \ar[r] \ar@{}[ru]|{push} & X\underset{B}{\cup}Y}
\end{displaymath}
Then $(X\underset{B}{\cup}Y,Y)$ is a relative CW($A$) (resp. a generalized relative CW($A$)).

Moreover, if $(X,B)$ has a CW($A$)-stucture of dimension $n\in\mathbb{N}_0$ (resp. a CW($A$)-structure with a finite number of layers) then $(X\underset{B}{\cup}Y,Y)$ can also be given a CW($A$)-stucture of dimension $n$ (resp. a CW($A$)-structure with a finite number of layers).
\end{lemma}

\begin{theo} \label{cwaconacwb1} 
Let $A$ be a CW($B$) of finite dimension and let $X$ be a generalized CW($A$). Then $X$ is a generalized CW($B$). In particular, if $A$ is a CW-complex of finite dimension then $X$ is a generalized CW-complex.
\end{theo}

\begin{proof}
Let $$\ast= X^0\rightarrow X^1 \rightarrow \ldots \rightarrow X^n \rightarrow \ldots$$ be a generalized CW($A$) structure on $X$.
Then, for each $n\in\mathbb{N}$ we have a pushout diagram
\begin{displaymath}
\xymatrix@C=30pt{C_n=\underset{\alpha\in J}{\bigvee}\Sigma^{n_\alpha-1}A \ar[r]^(.6){\underset{\alpha\in J}{+}g_\alpha} \ar[d]_i \ar@{}[rd]|{push} & X^{n-1} \ar[d] \\ D_n=(\underset{\alpha\in J_0}{\bigvee}A)\lor(\underset{\alpha\in J}{\bigvee}\cn\Sigma^{n_\alpha-1}A) \ar[r]_(.75){\underset{\alpha\in J}{+}f_\alpha} & X^n}
\end{displaymath}
where $n_\alpha\in\mathbb{N}$ for all $\alpha\in J$.

We have that $(D_n,C_n)$ is a relative CW($B$) by \ref{conoysuspcwa}, and it has finite dimension since $A$ does. So, by \ref{cwrelpaste}, $(X^n,X^{n-1})$ is a relative CW($B$) of finite dimension. Then, for each $n\in\mathbb{N}$, there exist spaces $Y^j_n$ for $0\leq j \leq m_n$, with $m_n\in\mathbb{N}$ such that $Y^j_n$ is obtained from $Y^{j-1}_n$ by attaching cells of type $B$ of dimension $j$ and $Y^{-1}_n=X^{n-1}$, $Y^{m_n}_n=X^n$. Thus, there exists a diagram
$$\ast= X^0=Y^{-1}_1\rightarrow Y^{0}_1 \rightarrow Y^{1}_1 \rightarrow \ldots \rightarrow Y^{m_1}_1 = X^1=Y^{-1}_2 \rightarrow \ldots \rightarrow Y^{m_2}_2 = X^2=Y^{-1}_3 \rightarrow \ldots$$
where each space is obtained from the previous one by attaching cells of type $B$. It is clear that $X$, the colimit of this diagram, is a generalized CW($B$).
\end{proof}

In the following example we exhibit a space $X$ which is not a CW-complex but is a CW($A$), with $A$ a CW-complex.

\begin{ex}
Let $A=[0;1]\cup\{2\}$, with $0$ as the base point. We build $X$ as follows. We attach two 0-cells to get $A\lor A$. We will denote the points in $A\lor A$ as $(a,j)$, where $a\in A$ and $j=1,2$. We define now, for each $n\in\mathbb{N}$, maps $g_n:A\rightarrow A\lor A$ in the following way. We set $g_n(a)=(a,1)$ if $a\in [0;1]$ and $g_n(2)=(1/n,2)$. We attach $1$-cells of type $A$ by means of the maps $g_n$. By a similar argument as the one in \ref{remcwavsgen}, the space $X$ obtained in this way is not a CW-complex.
\end{ex}

If $A$ is a finite dimensional CW-complex and $X$ is a generalized CW($A$), the previous theorem says that $X$ is a generalized CW-complex, and so it has the homotopy type of a CW-complex. The following result asserts that the last statement is also true for any CW-complex $A$.

\begin{prop} \label{homtype}
If $A$ is a CW-complex and $X$ is a generalized CW($A$) then $X$ has the homotopy type of a CW-complex.
\end{prop}

\begin{proof}
Let
\begin{displaymath}
\ast\subseteq X^1 \subseteq X^2 \subseteq \ldots \subseteq X^n\subseteq \ldots
\end{displaymath}
be a generalized CW($A$) structure on $X$. We may suppose that all the $0$-cells are attached in the first step, that is, $$X^1=\underset{\beta}{\bigvee}A\lor\underset{\alpha}{\bigvee}\Sigma^{n_\alpha}A$$
with $n_\alpha\in\mathbb{N}$. It is clear that $X^1$ is a CW complex.

We will construct inductively a sequence of CW-complexes $Y_n$ for $n\in\mathbb{N}$ with $Y_{n-1}\subseteq Y_n$ subcomplex and homotopy equivalences $\phi_n:X^n\rightarrow Y_n$ such that $\phi_n|_{X^{n-1}}=\phi_{n-1}$.

We take $Y_1=X^1$ and $\phi_1$ the identity map. Suppose we have already constructed $Y_1,\ldots,Y_k$ and $\phi_1,\ldots,\phi_k$ satisfying the conditions mentioned above. We consider the following pushout diagram.

\begin{displaymath}
\xymatrix{\underset{\alpha}{\bigvee}\Sigma^{n_\alpha-1}A \ar[r]^(.55){\underset{\alpha}{+}g_\alpha} \ar[d]_{\underset{\alpha}{\bigvee}i} \ar@{}[rd]|{push} & X^k \ar[d]^{i_k} \ar[r]^{\phi_k} \ar@{}[rd]|{push} & Y_k \ar[d]^{\gamma'_k} \\ \underset{\alpha}{\bigvee}\cn\Sigma^{n_\alpha-1}A \ar[r]_(.55){\underset{\alpha}{+}f_\alpha} &  X^{k+1} \ar[r]_\beta & Y'_{k+1} }
\end{displaymath}

Note that $\beta$ is a homotopy equivalence since $i_k$ is a closed cofibration and $\phi_k$ is a homotopy equivalence.

We deform $\phi_k\circ(\underset{\alpha}{+}g_\alpha)$ to a cellular map $\psi$ and we define $Y_{k+1}$ as the pushout

\begin{displaymath}
\xymatrix{\underset{\alpha}{\bigvee}\Sigma^{n_\alpha-1}A \ar[r]^(.6){\psi} \ar[d]_{\underset{\alpha}{\bigvee}i} \ar@{}[rd]|{push} & Y_k \ar[d]^{\gamma_k} \\ \underset{\alpha}{\bigvee}\cn\Sigma^{n_\alpha-1}A \ar[r] & Y_{k+1} }
\end{displaymath}

There exists a homotopy equivalence $k:Y'_{k+1}\rightarrow Y_{k+1}$ with $k|_{Y_k}=\id$. Let $i_k:X^k\rightarrow X^{k+1}$ be the inclusion. Then $k\beta i_k=k\gamma_k'\phi_k$ and $k\gamma_k'=\gamma_k$ is the inclusion. Let $\phi_{k+1}=k\beta$. Then, $\phi_{k+1}$ is a homotopy equivalence and $\phi_{k+1}|_{X^k}=\phi_k$.

We take $Y$ to be the colimit of the $Y_n$'s. Then $Y$ is a CW-complex. As the inclusions $i_k$, $\gamma_k$ are closed cofibrations, by proposition A.5.11 of \cite{FP}, it follows that $X$ is homotopy equivalent to $Y$.
\end{proof}

We prove now a variation of theorem \ref{cwaconacwb1}.

\begin{theo} \label{cwaconacwb2}
Let $A$ be a generalized CW($B$) with $B$ compact, and let $X$ be a generalized CW($A$). If $A$ and $B$ are T1 then $X$ is a generalized CW($B$).
\end{theo}

\begin{proof}
Let
$$\ast= X^0\rightarrow X^1 \rightarrow \ldots \rightarrow X^n \rightarrow \ldots$$
be a generalized CW($A$)-structure on $X$. Let $C_n$, $D_n$ be as in the proof of \ref{cwaconacwb1}.

We have that $(D_n,C_n)$ is a relative CW($B$) by \ref{conoysuspcwa}. By \ref{cwrelpaste}, $(X^n,X^{n-1})$ is 
also a relative CW($B$), but it need not be finite dimensional,
so we can not continue with the same argument as in the proof of \ref{cwaconacwb1}. But 
using the compactness of $B$, we will show that the cells of type $B$ may be attached in a 
certain order to obtain spaces $Z^n$ for $n\in\mathbb{N}$ such that $X$ is the colimit of the $Z^n$'s.

Let $J$ denote the set of all cells of type $B$ belonging to some of the relative 
CW($B$)'s $(X^n,X^{n-1})$ for $n\in\mathbb{N}$. We associate an ordered 
pair $(a,b)\in(\mathbb{N}_0)^2$ to each cell in $J$ in the following way. Note that each cell of type $B$ is 
included in exactly one cell of type $A$. The number $a$ will be the smallest number of layer in 
which that $A$-cell lies. In a similar way, if we regard that $A$-cell 
as a relative CW($B$) $(\cn\Sigma^{n-1}A,\Sigma^{n-1}A)$ (or more precisely, the image of this by the characteristic map), 
we set $b$ to be the smallest number of layer (in $(\cn\Sigma^{n-1}A,\Sigma^{n-1}A)$) in which the $B$-cell lies. 
If $e$ is the cell, we denote $\varphi(e)=(a,b)$.

We will consider in $(\mathbb{N}_0)^2$ the lexicographical order with the first coordinate greater than the second one.

Now we set the order in which the $B$-cells are attached. Let $J_1$ be the set of all the cells whose attaching 
map is the constant. We define inductively $J_n$ for $n\in\mathbb{N}$ to be the set of all the $B$-cells 
whose attaching map has image contained in the union of all the cells in $J_{n-1}$. Clearly $J_{n-1}\subseteq J_n$. 
We wish to attach first the cells of $J_1$, then those of $J_2-J_1$, etc. This can be done because of the 
construction of the $J_n$. We must verify that there are no cells missing, i.e., 
that $J=\underset{n\in\mathbb{N}}{\bigcup}J_n$.

Suppose there exists one cell in $J$, which we call $e_1$, which is not in any of the $J_n$. 
The image of its attaching map, denoted $K$,  is compact, since $B$ is compact and therefore
 it meets only a finite number of interiors of $A$-cells. For each of these cells $e_A$ we consider the 
relative CW($B$) $(\overline{e_A},\overline{e_A} - \overset{\circ}{e_A})$, where $e_A$ is the cell of type $A$. 

Then $K\cap\overline{e_A}$ is closed in $K$ and hence compact, so it meets only a finite number of interiors of 
$B$-cells of the relative CW($B$) $(\overline{e_A},\overline{e_A}-\overset{\circ}{e_A})$. 


Thus $K$ meets only a finite number of interiors of $B$-cells in $J$. 

This implies that $K$, which is the image of the attaching map of $e_1$, meets 
the interior of some cell $e_2$ which does not belong to any of the $J_n$, because of the finiteness condition. 

Recall that $e_2$ is an immediate face of $e_1$, which easily implies that $\varphi(e_2)<\varphi(e_1)$. 

Applying the same argument inductively we get a sequence of cells $(e_n)_{n\in\mathbb{N}}$ such 
that $\varphi(e_{n+1})<\varphi(e_n)$ for all $n$. 

But this induces an infinite decreasing sequence for the lexicographical order, which is impossible.
Hence, $J=\underset{n\in\mathbb{N}}{\bigcup}J_n$.

Let $Z^n=\underset{e\in J_n}{\bigcup}e$. It is clear that $(Z^n,Z^{n-1})$ is a relative CW($B$). 

Since colimits commute, we prove that $X=\col Z^n$ is a generalized CW($B$)-complex.  

\end{proof}

\section{The descriptive approach}


We will investigate now the descriptive approach and compare it with the constructive approach introduced
 in the previous section. We shall prove that in many cases a constructive CW($A$)-complex is the same as a descriptive one.

As before, let $A$ be a fixed pointed topological space.

\begin{definition}
Let $X$ be  a pointed topological space (with base point $x_0$). A \emph{cellular complex structure of type $A$ on $X$} is a collection $\mathcal{K}=\{e^n_\alpha:n\in\mathbb{N}_0,\alpha\in J_n\}$ of subsets of $X$, which are called the \emph{cells (of type $A$)}, such that $x_0\in e^n_\alpha$ for all $n$ and $\alpha$, and satisfying conditions (1), (2) and (3) below.

Let $\mathcal{K}^n=\{e^r_\alpha,r\leq n, \alpha\in J_r\}$ for $n\in\mathbb{N}_0$, $\mathcal{K}^{-1}=\{\{x_0\}\}$. $\mathcal{K}^n$ is called the \emph{$n$-skeleton of $\mathcal{K}$}. Let $|\mathcal{K}^n|=\underset{\substack{r\leq n \\ \alpha\in J_r} } {\bigcup}e^r_\alpha$, $|\mathcal{K}^n|\subseteq X$ a subspace.

We call $\overset{\bullet}{e^n_\alpha}= e^n_\alpha \cap |\mathcal{K}^{n-1}|$ the \emph{boundary of the cell $e^n_\alpha$} and $\overset{\circ}{e^n_\alpha}= e^n_\alpha - \overset{\bullet}{e^n_\alpha}$ the \emph{interior of the cell $e^n_\alpha$}.

The collection $\mathcal{K}$ must satisfy the following properties.
\begin{enumerate}
\item[(1)] $X=\underset{n,\alpha}{\bigcup}e^n_\alpha=|\mathcal{K}|$
\item[(2)] $\overset{\circ}{e^n_\alpha}\cap \overset{\circ}{e^m_\beta}\neq \varnothing \Rightarrow m=n, \alpha=\beta$
\item[(3)] For every cell $e^n_\alpha$ with $n\geq 1$ there exists a continuous map $$f^n_\alpha:(\cn\Sigma^{n-1}A,\Sigma^{n-1}A,a_0)\rightarrow (e^n_\alpha,\overset{\bullet}{e^n_\alpha},x_0)$$ such that $f^n_\alpha$ is surjective and $f^n_\alpha:\cn\Sigma^{n-1}A-\Sigma^{n-1}A \rightarrow \overset{\circ}{e^n_\alpha}$ is a homeomorphism. For $n=0$, there is a homeomorphism $f^0_\alpha:(A,a_0)\rightarrow (e^0_\alpha,x_0)$.
\end{enumerate}
The \emph{dimension of $\mathcal{K}$} is defined as $\dim \mathcal{K}= \sup\{n:J_n\neq\varnothing\}$.
\end{definition}

\begin{definition}
Let $\mathcal{K}$ be a cellular complex structure of type $A$ in a topological space $X$. We say that \emph{$\mathcal{K}$ is a cellular CW-complex with base $A$} if it satisfies the following conditions.
\begin{enumerate}
\item[(C)] Every compact subspace of $X$ intersects only a finite number of interiors of cells.
\item[(W)] $X$ has the weak (final) topology with respect to the cells.
\end{enumerate}
In this case we will say that $X$ is a \emph{descriptive CW($A$)}.
\end{definition}

We study now the relationship between both approaches. 

\begin{theo}
Let $A$ be a T1 space. If $X$ is a constructive CW($A$), then it is a descriptive CW($A$).
\end{theo}

\begin{proof}
Let $\mathcal{K}=\{e^n_\alpha\}_{n,\alpha} \cup\{\{x_0\}\}$. It is not difficult to verify that $\mathcal{K}$ defines a cellular complex structure on $X$.

It remains to prove that it satisfies conditions (C) and (W). Note that condition (C) follows from \ref{compintersfinitosint}, while (W) follows from \ref{remweaktopo}.
%
\end{proof}

Note that the hypothesis of T1 on $A$ is necessary. For example, take $A=\{0,1\}$ with the indiscrete topology and 0 as base point. Let $X=\underset{j\in\mathbb{N}}{\bigvee}A$. The space $X$ also has the indiscrete topology and it is a constructive CW($A$). If it were a descriptive CW($A$), it could only have cells of dimension 0 since $X$ is countable. But $X$ is not finite, then it must have infinite many cells, but it is a compact space. This implies that (C) does not hold, thus $X$ is not a descriptive CW($A$).

\begin{theo} \label{desc_impl_const}
Let $A$ be a compact space and let $X$ be a descriptive CW($A$). If $X$ is Hausdorff then it is a constructive CW($A$).
\end{theo}

\begin{proof}
We will prove that $|\mathcal{K}^n|$ can be obtained from $|\mathcal{K}^{n-1}|$ by attaching $A$-$n$-cells. For $n=0$ this is clear since we have a homeomorphism $\underset{\alpha\in J_0}{\bigvee}f^0_\alpha:\underset{\alpha\in J_0}{\bigvee}A \rightarrow |\mathcal{K}^0|$.

For any $n\in\mathbb{N}$, there is a pushout

\begin{displaymath}
\xymatrix@C=45pt{\underset{\alpha\in J_n}{\bigvee}\Sigma^{n-1}A \ar[r]^(.6){\underset{\alpha\in J_n}{+}f^n_\alpha|_{\Sigma^{n-1}A}} \ar[d]_i \ar@{}[rd]|{push} & |\mathcal{K}^{n-1}| \ar[d] \\ \underset{\alpha\in J_n}{\bigvee}\cn\Sigma^{n-1}A \ar[r]_(.6){\underset{\alpha\in J_n}{+}f^n_\alpha} & |\mathcal{K}^n|}
\end{displaymath}

The topology of $|\mathcal{K}^n|$ coincides with the pushout topology since X is 
hausdorff and A is compact.
\end{proof}

It is interesting to see that \ref{desc_impl_const} is not true if $X$ is not Hausdorff, even in the case $A$ is compact and Hausdorff. For example, take $A=S^0$ with the usual topology, and $X=[-1;1]$ with the following topology. The proper open sets are $[-1;1)$, $(-1;1]$ and the subsets $U\subseteq (-1;1)$ which are open in $(-1;1)$ with the usual topology. It is easy to see that $X$ is a descriptive CW($A$). We denote $D^1=[-1;1]$ with the usual topology. Take $e^0=\{-1;1\}$, $e^1=X$. Let $f^0:A\rightarrow \{-1;1\}$ and $f^1:CA=D^1\rightarrow e^1$ be the identity maps on the underlying sets. Both maps are continuous and surjective. The maps $f^0$ and $f^1|_{\overset{\circ}{D^1}}:\overset{\circ}{D^1}\rightarrow \overset{\circ}{e^1}$ are homeomorphisms. So conditions (1), (2) and (3) of the definition of cellular complex are satisfied. Condition (C) is obvious, and (W) follows from the fact that $e^1=X$. So $X$ is a descriptive CW($A$). But it is not a constructive CW($A$) because it is not Hausdorff.

\bigskip

In a similar way one can define the notion of \emph{descriptive generalized CW($A$)-complex}. The relationship between the constructive and descriptive approachs of generalized CW($A$)-complexes is analogous to the previous one.

\bigskip

\section{Changing cores}

Suppose we have two spaces $A$ and $B$ and maps $\alpha:A\rightarrow B$ and $\beta:B\rightarrow A$. Let $X$ be a CW($A$). We want to construct a CW($B$) out of $X$, using the maps $\alpha$ and $\beta$.

We shall consider two special cases. First, we consider the case $\beta\alpha=\id_A$, that is, $A$ is a retract of $B$. In this case, we construct a CW($B$) $Y$ such that $X$ is a retract of $Y$. 

\medskip

We denote $g_\gamma^n$, $f_\gamma^n$ the adjunction and characteristic maps of the $A$-$n$-cells ($\gamma \in J_n$). Let $Y^0=\underset{\gamma\in J_0}{\bigvee}B$ and let $\varphi_0:X^0\rightarrow Y^0$ be the map $\lor\alpha$ and let $\psi_0:Y^0\rightarrow X^0$ be the map $\lor\beta$. Clearly $\psi_0\varphi_0=\id_{X^0}$.

By induction suppose we have constructed $Y^{n-1}$ and maps $\varphi_{n-1}:X^{n-1}\rightarrow Y^{n-1}$ and $\psi_{n-1}:Y^{n-1}\rightarrow X^{n-1}$ such that $\psi_{n-1}\varphi_{n-1}=\id_{X^{n-1}}$ and such that $\varphi_{k}$, $\psi_{k}$ extend $\varphi_{k-1}$, $\psi_{k-1}$ for all $k\leq n-1$. We define $Y^n$ by the following pushout.

\begin{displaymath}
\xymatrix@C=65pt{\underset{\gamma\in J_n}{\bigvee}\Sigma^{n-1}B \ar[r]^(.6){\varphi_{n-1}(\underset{\gamma\in J_n}{+}g^n_\gamma\Sigma^{n-1}\beta)} \ar[d]_{\lor i} \ar@{}[rd]|{push} & Y^{n-1} \ar[d]^j \\ \underset{\gamma\in J_n}{\bigvee}\cn\Sigma^{n-1}B \ar[r]_(.6){\underset{\gamma\in J_n}{+}h^n_\gamma} & Y^n}
\end{displaymath}

Since $$\begin{array}{ll}
(\underset{\gamma\in J_n}{+}f^n_\gamma\cn\Sigma^{n-1}\beta)(\lor i) & =\underset{\gamma\in J_n}{+}(f^n_\gamma\cn\Sigma^{n-1}\beta i) = \underset{\gamma\in J_n}{+}(f^n_\gamma i\Sigma^{n-1}\beta)= \underset{\gamma\in J_n}{+}(\incl g^n_\gamma \Sigma^{n-1}\beta) = \\ & = \incl \psi_{n-1} \underset{\gamma\in J_n}{+}(\varphi_{n-1} g^n_\gamma \Sigma^{n-1}\beta)\end{array}$$

there exists a map $\psi_n:Y^n\rightarrow X^n$ extending $\psi_{n-1}$ such that $\psi_n \underset{\gamma\in J_n}{+}h^n_\gamma = \underset{\gamma\in J_n}{+}(f^n_\gamma \cn\Sigma^{n-1}\beta)$ and $\psi_n j = \incl\psi_{n-1}$.

On the other hand we have the following commutative diagram

\begin{displaymath}
\xymatrix@C=20pt@R=20pt{\underset{\gamma\in J_n}{\bigvee}\Sigma^{n-1}A \ar[rr]^(.5){\underset{\gamma\in J_n}{+}g^n_\gamma} \ar[dd]_{\lor i} \ar[rd]^{\lor\Sigma^{n-1}\alpha} & & X^{n-1} \ar[dd]^(.65){\incl} \ar[rd]^{\varphi_{n-1}} \\ 
& \underset{\gamma\in J_n}{\bigvee}\Sigma^{n-1}B \ar[rr]^(.55){\varphi_{n-1}(\underset{\gamma\in J_n}{+}g^n_\gamma\Sigma^{n-1}\beta)} \ar[dd]_(.35){\lor i} & & Y^{n-1} \ar[dd]^j \\
\underset{\gamma\in J_n}{\bigvee}\cn\Sigma^{n-1}A \ar[rr]_(.7){\underset{\gamma\in J_n}{+}f^n_\gamma} \ar[rd]^{\lor\cn\Sigma^{n-1}\alpha} & & X^n \ar@{.>}[rd]^{\varphi_n} & \\
& \underset{\gamma\in J_n}{\bigvee}\cn\Sigma^{n-1}B \ar[rr]_(.5){\underset{\gamma\in J_n}{+}h^n_\gamma} & & Y^n
}
\end{displaymath}
where the front and back faces are pushouts. Then the dotted arrow exists and we have $\varphi_n=j\varphi_{n-1}+(\underset{\gamma\in J_n}{+}h^n_\gamma\cn\Sigma^{n-1}\alpha)$. Also, $\psi_n\varphi_n=\id_{X^n}$, since $$\begin{array}{ll}\psi_n\varphi_n & =\psi_nj\varphi_{n-1}+(\underset{\gamma\in J_n}{+}\psi_nh^n_\gamma\cn\Sigma^{n-1}\alpha)=\incl\psi_{n-1}\varphi_{n-1} + (\underset{\gamma\in J_n}{+}f^n_\gamma\cn\Sigma^{n-1}\beta\cn\Sigma^{n-1}\alpha)= \\ & = \incl + (\underset{\gamma\in J_n}{+}f^n_\gamma) = \id_{X^n}\end{array}$$

Let $Y=\col Y^n$. Then there exist maps $\varphi:X\rightarrow Y$ and $\psi:Y\rightarrow X$ induced by the $\psi_n$'s and $\varphi_n$'s and they satisfy $\psi\varphi=\id_X$. So, $X$ is a retract of $Y$.

\bigskip

The second special case we consider is the following. Suppose $A$ and $B$ have the same homotopy type, that is, there exists a homotopy equivalence $\beta:B\rightarrow A$ with homotopy inverse $\alpha$. Suppose, in addition, that the base points of $A$ and $B$ are closed. Let $X$ be a CW($A$). We will construct a CW($B$) which is homotopy equivalent to $X$.

Again we take $Y^0=\underset{\gamma\in J_0}{\bigvee}B$. Let $\varphi_0:X^0\rightarrow Y^0$ be the map $\lor\alpha$. So, $\varphi_0$ is a homotopy equivalence. 

Now, let $n\in\mathbb{N}$ and suppose we have constructed $Y^{n-1}$ and a homotopy equivalence $\varphi_{n-1}:X^{n-1}\rightarrow Y^{n-1}$. We define $Y^n$ as in the first case. Consider the commutative diagrams
\begin{displaymath}
\xymatrix@C=15pt@R=15pt{\underset{\gamma\in J_n}{\bigvee}\Sigma^{n-1}B \ar[rr]^(.5){\underset{\gamma\in J_n}{+}g^n_\gamma\Sigma^{n-1}\beta} \ar[dd]_{i_B} \ar[rd]^{\id} & & X^{n-1} \ar[dd]^(.65){\incl} \ar[rd]^{\varphi_{n-1}} \\
& \underset{\gamma\in J_n}{\bigvee}\Sigma^{n-1}B \ar[rr]^(.55){\varphi_{n-1}(\underset{\gamma\in J_n}{+}g^n_\gamma\Sigma^{n-1}\beta)} \ar[dd]_(.35){i_B} & & Y^{n-1} \ar[dd]^j \\
\underset{\gamma\in J_n}{\bigvee}\cn\Sigma^{n-1}B \ar[rr] \ar[rd]^{\id} & & X^{n-1}\cup e^n_B \ar@{.>}[rd]^{p_1} & \\
& \underset{\gamma\in J_n}{\bigvee}\cn\Sigma^{n-1}B \ar[rr]_(.5){\underset{\gamma\in J_n}{+}h^n_\gamma} & & Y^n
}
\end{displaymath}
\begin{displaymath}
\xymatrix@C=15pt@R=15pt{\underset{\gamma\in J_n}{\bigvee}\Sigma^{n-1}B \ar[rr]^(.5){\underset{\gamma\in J_n}{+}g^n_\gamma\Sigma^{n-1}\beta} \ar[dd]_{i_B} \ar[rd]^{\lor\Sigma^{n-1}\beta} & & X^{n-1} \ar[dd]^(.65){\incl} \ar[rd]^{\id} \\
& \underset{\gamma\in J_n}{\bigvee}\Sigma^{n-1}A \ar[rr]^(.4){\underset{\gamma\in J_n}{+}g^n_\gamma} \ar[dd]_(.35){i_A} & & X^{n-1} \ar[dd]^\incl \\
\underset{\gamma\in J_n}{\bigvee}\cn\Sigma^{n-1}B \ar[rr] \ar[rd]_{\lor\cn\Sigma^{n-1}\beta} & & X^{n-1}\cup e^n_B \ar@{.>}[rd]^{p_2} & \\
& \underset{\gamma\in J_n}{\bigvee}\cn\Sigma^{n-1}A \ar[rr]_(.5){\underset{\gamma\in J_n}{+}f^n_\gamma} & & X^n
}
\end{displaymath}

Since the front and rear faces of both cubical diagrams are pushouts, the dotted arrows $p_1$ and $p_2$ exist. Now $\varphi_{n-1}$, $\lor\Sigma^{n-1}\beta$ and $\lor\cn\Sigma^{n-1}\beta$ are homotopy equivalences and $i_A$ and $i_B$ are closed cofibrations. Then, by proposition 7.5.7 of \cite{Br}, $p_1$ and $p_2$ are homotopy equivalences. 
We have the following commutative diagram.

\begin{displaymath}
\xymatrix{ Y^{n-1} \ar[d]^{i} & X^{n-1} \ar[l]_{\varphi_{n-1}} \ar[r]^\id \ar[d]^{j} & X^{n-1} \ar[d]^{k} \\ Y^n & X^{n-1}\cup e^n_B \ar[l]_{p_1} \ar[r]^{p_2} & X^n }
\end{displaymath}
where $i$, $j$ and $k$ are the inclusions. Let $p_2^{-1}$ be a homotopy inverse of $p_2$. Then $p_1p_2^{-1}k= p_1p_2^{-1}p_2j \simeq p_1j =i\varphi_{n-1}$. Since $k:X^{n-1}\rightarrow X^n$ is a cofibration, $\varphi_{n-1}$ extends to some $\varphi_{n}:X^n\rightarrow Y^n$ and $\varphi_n$ is homotopic to $p_1p_2^{-1}$, and thus, it is a homotopy equivalence.

Again, we take $Y=\col Y^n$. Then the maps $\varphi_n$ for $n\in\mathbb{N}$ induce a map $\varphi:X\rightarrow Y$ which is a homotopy equivalence by proposition A.5.11 of \cite{FP}.

We summarize the previous results in the following theorem.

\begin{theo}
Let $A$ and $B$ be pointed topological spaces. Let $X$ be a CW($A$), and let $\alpha:A\rightarrow B$ and $\beta:B\rightarrow A$ be continuous maps.
\begin{enumerate}
\item[i.] If $\beta\alpha=\id_A$, then there exists a CW($B$) $Y$ and maps $\varphi:X\rightarrow Y$ and $\psi:Y\rightarrow X$ such that $\psi\varphi=\id_X$.
\item[ii.] Suppose $A$ and $B$ have closed base points. If $\beta$ is a homotopy equivalence, then there exists a CW($B$) $Y$ and a homotopy equivalence $\varphi:X\rightarrow Y$.
\item[iii.] Suppose $A$ and $B$ have closed base points. If $\beta\alpha=\id_A$ and $\alpha\beta\simeq\id_A$ then there exists a CW($B$) $Y$ and maps $\varphi:X\rightarrow Y$ and $\psi:Y\rightarrow X$ such that $\psi\varphi=\id_X$ and $\varphi\psi\simeq\id_Y$.
\end{enumerate}
\end{theo}

Note that item (iii) follows by a similiar argument.

\medskip

The previous theorem has an easy but interesting corollary.

\begin{coro}
Let $A$ be a contractible space (with closed base point) and let $X$ be a CW($A$). Then $X$ is contractible.
\end{coro}

This corollary also follows from a result analogous to Whitehead Theorem which we prove in the next section.

\section{Homotopy theory of CW(A)-complexes}

In this section we start to develop the homotopy theory of CW($A$)-complexes. The main result of this section is theorem \ref{Wh_theo} which generalizes the famous Whitehead Theorem.

\bigskip

Let $X$ be a (pointed) topological space and let $r\in\mathbb{N}_0$. Recall that the sets $\pi^A_r(X)$ are defined by $\pi^A_r(X)=[\Sigma^r A,X]$, the homotopy classes of maps from $\Sigma^r A$ to $X$. It is well known that these are groups for $r\geq 1$ and Abelian for $r\geq 2$.

Similarly, for $B\subseteq X$ one defines $\pi^A_r(X,B)=[(\cn\Sigma^{r-1} A,\Sigma^{r-1} A),(X,B)]$ for $r\in\mathbb{N}$, which are groups for $r\geq 2$ and Abelian for $r\geq 3$.

Note that $\pi^{S^0}_r(X)=\pi_r(X)$ and $\pi^{S^n}_r(X)=\pi_{r+n}(X)$. Note also that $\pi^A_r(X)$ are trivial if $A$ is contractible.

\begin{definition}
Let $(X,B)$ be a pointed topological pair. The pair $(X,B)$ is called \emph{$A$-0-connected} if for any given continuous function $f:A\rightarrow X$ there exists a map $g:A\rightarrow B$ such that $ig\simeq f$, where $i:B\rightarrow X$ is the inclusion.
\begin{displaymath}
\xymatrix{\ast \ar[r] \ar[d] \ar@{}[rd]|(.65){\simeq} & B \ar[d]^i \\ A \ar[r]_f \ar[ru]^g & X}
\end{displaymath}
\end{definition}

\begin{definition}
Let $n\in\mathbb{N}$. The pointed topological pair $(X,B)$ is called \emph{$A$-n-connected} if it is $A$-0-connected and $\pi_r^A(X,B)=0$ for $1\leq r \leq n$.
\end{definition}

\begin{definition}
Let $f:X\rightarrow Y$ be a continuous map, and let $A$ be a topological space. The map $f$ is called an \emph{$A$-0-equivalence} if for any given continuous function $g:A\rightarrow Y$, there exists a map $h:A\rightarrow X$ such that $fh\simeq g$.
\begin{displaymath}
\xymatrix{\ast \ar[r] \ar[d] \ar@{}[rd]|(.65){\simeq} & X \ar[d]^f \\ A \ar[r]_g \ar[ru]^h & Y}
\end{displaymath}
Given $n\in\mathbb{N}$, the map $f$ is called an \emph{$A$-n-equivalence} if it induces isomorphisms $f_*:\pi_r^A(X,x_0)\rightarrow \pi_r^A(Y,f(x_0))$ for $0\leq r < n$ and an epimorphism for $r=n$.
\newline
Also, $f$ is called an \emph{$A$-weak equivalence} if it is an $A$-$n$-equivalence for all $n\in\mathbb{N}$.
\end{definition}

\begin{rem}
Let $f:X\rightarrow Y$ be map and let $n\in\mathbb{N}$. We denote by $Z_f$ the mapping cylinder of $f$. Then $f$ is an \emph{$A$-$n$-equivalence} if and only if the topological pair $(Z_f,X)$ is $A$-$n$-connected.
\end{rem}

\begin{lemma}\label{relhomotlemma}
Let $X$, $S$, $B$ be pointed topological spaces, $S\subseteq X$ a subspace, $x_0\in S$ and $b_0\in B$ the base points. Let $f:(\cn B,B)\rightarrow(X,S)$ be a continuous map. Then the following are equivalent.
\begin{enumerate}
\item[i)] There exists a base point preserving homotopy $H:(\cn B\times I,B\times I)\rightarrow(X,S)$ such that $Hi_0=f$, $Hi_1(x)=x_0$ $\forall x\in\cn B$.
\item[ii)] There exists a (base point preserving) homotopy $G:\cn B\times I\rightarrow X$, relative to $B$, such that $Gi_0=f$, $Gi_1(\cn B)\subseteq S$.
\item[iii)] There exists a (base point preserving) homotopy $G:\cn B\times I\rightarrow X$, such that $Gi_0=f$, $Gi_1(\cn B)\subseteq S$.
\end{enumerate}
\end{lemma}

\begin{proof}
$i)\Rightarrow ii)$ Define $G$ as follows.
\begin{displaymath}
G([x,s],t)= \left\{ \begin{array}{ll} H([x,\frac{2s}{2-t}],t) & \mathrm{if}\ 0\leq s\leq 1-\frac{t}{2} \\ H([x,1],2-2s) & \mathrm{if}\ 1-\frac{t}{2}\leq s\leq 1 \end{array} \right.
\end{displaymath}
It is clear that $G$ is well defined and continuous. Note that
\begin{displaymath}
\begin{array}{ll}
Gi_0([x,s])=H([x,\frac{2s}{2}],0)=H([x,s],0)=f(x,s) \\
Gi_1([x,s])=H([x,2s],1)=x_0 \in S & \mathrm{if}\ s \leq \frac{1}{2} \\
Gi_1([x,s])=H([x,1],2-2s) \in S & \mathrm{if}\ s \geq \frac{1}{2}
\end{array}
\end{displaymath}
since $H(B\times I)\subseteq S$.
\newline
$ii)\Rightarrow iii)$ Obvious.
\newline
$iii)\Rightarrow i)$ We define $H$ by
\begin{displaymath}
H([x,s],t) = \left\{ \begin{array}{ll} G([x,s],2t) & \mathrm{if}\ 0\leq t\leq \frac{1}{2} \\ Gi_1([x,s(2-2t)]) & \mathrm{if}\ \frac{1}{2}\leq t\leq 1 \end{array} \right.
\end{displaymath}
\end{proof}

\begin{lemma}
Let $X$, $Y$ be pointed topological spaces 
and let $f:X\rightarrow Y$ be an $A$-$n$-equivalence. Let $r\in\mathbb{N}$, $r\leq n$ and let $i_A:\Sigma^{r-1}A\rightarrow\cn\Sigma^{r-1}A$ be the inclusion. Suppose that $g:\Sigma^{r-1}A\rightarrow X$ and $h:\cn\Sigma^{r-1}A\rightarrow Y$ are continuous maps such that $hi_A=fg$. Then, there exists a continuous map $k:\cn\Sigma^{r-1}A\rightarrow X$ such that $ki_A=g$ and $fk\simeq h \ \mathrm{rel}\ \Sigma^{r-1}A$.
\begin{displaymath}
\xymatrix{\Sigma^{r-1}A \ar[r]^(.6){g} \ar[d]_{i_A} \ar@{}[rd]|(.65){\simeq}|(.35){\mid\parallel} & X \ar[d]^f \\ \cn\Sigma^{r-1}A \ar[r]_(.6){h} \ar[ru]^(.6){k} & Y}
\end{displaymath}
\end{lemma}

\begin{proof}
Consider the inclusions $i:X\rightarrow Z_f$ and $j:Y\rightarrow Z_f$. Let $r:Z_f\rightarrow Y$ be the usual retraction. Note that there is a homotopy commutative diagram
\begin{displaymath}
\xymatrix{\Sigma^{r-1}A \ar[r]^(.6){g} \ar[d]_{i_A} & X \ar[d]^i \\ \cn\Sigma^{r-1}A \ar[r]_(.6){jh} & Z_f}
\end{displaymath}
Let $H:\Sigma^{r-1}A\times I \rightarrow Z_f$ be the homotopy from $jhi_A$ to $ig$ defined by $H(a,t)=[g(a),t]$ for $a\in\Sigma^{r-1}A$, $t\in I$. Consider the commutative diagram of solid arrows
\begin{displaymath}
\xymatrix{\Sigma^{r-1}A \ar[r]^{i_0} \ar[d]_{i_A} & \Sigma^{r-1}A\times I \ar[d] \ar@/^3ex/[rdd]^H & \\ \cn\Sigma^{r-1}A \ar[r]_{i_0} \ar@/_3ex/[rrd]^{jh} & \cn\Sigma^{r-1}A \times I \ar@{.>}[rd]^{H'} & \\ & & Z_f}
\end{displaymath}
Since $i_A$ is a cofibration there exists a map $H'$ such that the whole diagram commutes, which induces a commutative diagram
\begin{displaymath}
\xymatrix{\Sigma^{r-1}A \ar[r]^(.6){g} \ar[d]_{i_A} & X \ar[d]^i \\ \cn\Sigma^{r-1}A \ar[r]_(.6){H'i_1} & Z_f}
\end{displaymath}
The pair $(Z_f,X)$ is $A$-$n$-connected, so by lemma \ref{relhomotlemma} there exists a continuous function $k:\cn\Sigma^{r-1}A\rightarrow X$ such that $ki_A=g$, $ik\simeq H'i_1$ rel $\Sigma^{r-1}A$.
Then $$fk=rik\simeq rH'i_1\simeq rH'i_0=rjh=h$$ 
Note that the homotopy is relative to $\Sigma^{r-1}A$, thus $fk\simeq h$ rel $\Sigma^{r-1}A$.
\end{proof}

\begin{theo} \label{levnequiv}
Let $f:X\rightarrow Y$ be an $A$-$n$-equivalence ($n=\infty$ is allowed) and let $(Z,B)$ be a relative CW($A$) which admits a CW($A$)-structure of dimension less than or equal to $n$. Let $g:B\rightarrow X$ and $h:Z\rightarrow Y$ be continuous functions such that $h|_B=fg$. Then there exists a continuous map $k:Z\rightarrow X$ such that $k|_B=g$ and $fk\simeq h$ rel $B$.
\begin{displaymath}
\xymatrix{B \ar[r]^g \ar[d]_i \ar@{}[rd]|(.65){\simeq}|(.35){\mid\parallel} & X \ar[d]^f \\ Z \ar[r]_h \ar@{.>}[ru]|{k} & Y}
\end{displaymath}
\end{theo}

\begin{proof}
Let 
$$\begin{array}{lcl}S & = & \{(Z',k',K')/B\subseteq Z'\subseteq Z\ A-\textrm{subcomplex },\  k':Z'\rightarrow Z \textrm{ with } k'|_{B}=g \textrm{ and } \\ & &  K':Z'\times I\rightarrow Y , K':fk'\simeq h|_{Z'} \textrm{ rel } B \} \end{array}$$
It is clear that $S\neq\varnothing$. We define a partial order in $S$ in the following way.
$$(Z',k',K')\leq(Z'',k'',K'') \ \textrm{if and only if}\ Z'\subseteq Z'', k''|_{Z'}=k' \ \mathrm{K''|_{Z'\times I}=K'}$$
It is clear that every chain has an upper bound since $Z$ has the weak topology.
Then, by Zorn's lemma, there exists a maximal element $(Z',k',K')$. We want to prove that $Z'=Z$. Suppose $Z'\neq Z$, then there exist some $A$-cells in $Z$ which are not in $Z'$. Choose $e$ an $A$-cell with minimum dimension. We want to extend the maps $k'$ and $K'$ to $Z'\cup e$. If $e$ is an $A$-0-cell this is easy to do since $f$ is an $A$-0-equivalence and all homotopies are relative to the base point. Suppose then that $\dim e\geq 1$. Let $\phi:(\cn\Sigma^{r-1}A, \Sigma^{r-1}A)\rightarrow (Z,Z')$ be the characteristic map of $e$, let $\psi=\phi|_{\Sigma^{r-1}A}$, and let $Z''=Z'\cup e$. We have the following diagram.
\begin{displaymath}
\xymatrix{\Sigma^{r-1}A \ar[r]^\psi \ar[d]_{i_A} & Z' \ar[r]^{k'} \ar[d]_{i_{Z'}} & X \ar[d]^f \\ \cn\Sigma^{r-1}A \ar[r]_\phi \ar@{}[ru]|{\mid\parallel} & Z'' \ar[r]_{h|_{Z''}} \ar@{}[ru]|{\simeq} & Y}
\end{displaymath}
Here, the homotopy of the right square is relative to $B$. Let $\alpha:I\rightarrow I$ be defined by $\alpha(t)=1-t$. Since $i_{Z'}$ is a cofibration we can extend $K'(\id\times\alpha)$ to some $H:Z''\times I\rightarrow Y$, and then we obtain a commutative diagram
\begin{displaymath}
\xymatrix{\Sigma^{r-1}A \ar[r]^\psi \ar[d]_{i_A} & Z' \ar[r]^{k'} \ar[d]_{i_{Z'}} & X \ar[d]^f \\ \cn\Sigma^{r-1}A \ar[r]_\phi \ar@{}[ru]|{\mid\parallel} & Z'' \ar[r]_{Hi_1} \ar@{}[ru]|{\mid\parallel} & Y}
\end{displaymath}
By the previous lemma, there exists $l:\cn\Sigma^{r-1}A\rightarrow X$ such that $li_A=k'\psi$ and $fl\simeq Hi_1\phi$ rel $\Sigma^{r-1}A$. Let $G$ denote this homotopy.
\newline
Now, since the left square is a pushout, there is a map $\gamma:Z''\rightarrow X'$ such that $\gamma\phi=l$, $\gamma i_{Z'}=k'$. So $\gamma$ extends $k'$. We want now to define a homotopy $K'':f\gamma\simeq h|_{Z''}$ extending $K'$. 
We consider $\cn\Sigma^{r-1}A\times [0,2]/\sim$ where we identify $(b,t)\sim(b,t')$ for $b\in\Sigma^{r-1}A$, $t,t'\in[1,2]$. There is a homeomorphism $\beta:\cn\Sigma^{r-1}A\times [0,2]/\sim \rightarrow \cn\Sigma^{r-1}A\times I$ defined by 
\begin{displaymath}
\beta([a,s],t)=\left\{\begin{array}{ll} ([a,s],\frac{t}{2-s}) & \textrm{if}\ 0\leq t\leq 1 \\ ([a,s],\frac{1-s}{2-s}\ t+\frac{s}{2-s}) & \textrm{if}\ 1\leq t\leq 2 \end{array} \right.
\end{displaymath}
We have the following commutative diagram.
\begin{displaymath}
\xymatrix{\Sigma^{r-1}A\times I \ar[r]^(.58){\psi\times\id_I} \ar[d]_{i_A\times \id_I} & Z'\times I \ar[d]^{i_{Z'}\times \id_I} \ar@/^3ex/[rdd]^{K'(\id\times\alpha)} \\ \cn\Sigma^{r-1}A\times I \ar[r]_(.58){\phi\times\id_I} \ar@{}[ru]|{\textrm{push}} \ar@/_3ex/[rrd]_{(H(\phi\times\id_I)+G(\id\times\alpha))\beta^{-1}} & Z''\times I \ar@{.>}[rd]_{\overset{\sim}{K}} \\ & & Y}
\end{displaymath}

Note that
\begin{displaymath}
\begin{array}{l}
(H(\phi\times\id_I)+G(\id\times\alpha))\beta^{-1}(i_A\times\id_I) = H(\phi\times\id_I)(i_A\times\id_I) = \\ = H(i_{Z'}\times\id_I)(\psi\times\id_I) = K' (\id\times\alpha)(\psi\times\id_I)
\end{array}
\end{displaymath}

Then, the map $\overset{\sim}{K}$ exists. We take $K''=\overset{\sim}{K}(\id\times\alpha)$.
\end{proof}

\begin{rem}
If $(Y,B)$ is a relative CW($A$) which is $A$-$n$-connected for all $n\in\mathbb{N}$ then $i:B\rightarrow Y$ is an $A$-$n$-equivalence for all $n\in\mathbb{N}$ and we have
\begin{displaymath}
\xymatrix{B \ar[r]^{\id_B} \ar[d]_i \ar@{}[rd]|(.65){\simeq}|(.35){\mid\parallel} & B \ar[d]^i \\ Y \ar[r]_{\id_Y} \ar@{.>}[ru]|{r} & Y}
\end{displaymath}
Thus $B$ is a strong deformation retract of $Y$. In particular, if $X$ is a CW($A$) with $\pi_n^A(X)=0$ for all $n\in\mathbb{N}_0$, then $X$ is contractible.
\end{rem}

The following proposition follows immediately from \ref{levnequiv}.

\begin{prop}
Let $f:Z\rightarrow Y$ be an $A$-$n$-equivalence ($n=\infty$ is allowed) and let $X$ be a CW($A$) which admits a CW($A$)-structure of dimension less than or equal to $n$. Then, the map $f_*:[X,Z]\rightarrow[X,Y]$ is surjective.
\end{prop}

Finally we obtain a generalization of Whitehead's theorem.

\begin{theo} \label{Wh_theo}
Let $X$, $Y$ be CW($A$)'s and $f:X\rightarrow Y$ a continuous map. Then $f$ is a homotopy equivalence if and only if it is an $A$-weak equivalence.
\end{theo}

\begin{proof}
Suppose $f$ is an $A$-weak equivalence. We consider $f_*:[Y,X]\rightarrow[Y,Y]$. By the previous proposition, $f_*$ is surjective, then there exists $g:Y\rightarrow X$ such that $fg\simeq\id_Y$. Then $g$ is also an $A$-weak equivalence, so applying the above argument, there exists an $h:X\rightarrow Y$ such that $gh\simeq\id_X$. Then $f\simeq fgh\simeq h$, and so, $gf\simeq gh \simeq \id_X$. Thus $f$ is a homotopy equivalence.
\end{proof}

We finish with some results concerning the connectedness of CW($A$)-complexes.

\begin{lemma}
Let $A$ be an $l$-connected CW-complex, let $B$ be a topological space, and suppose $X$ is obtained from $B$ by attaching a 1-cell of type $A$. Then $(X,B)$ is $(l+1)$-connected.
\end{lemma}

\begin{proof}
Let $g$ be the attaching map of the cell and $f$ its characteristic map. Since $A$ is an $l$-connected CW-complex, $(\cn A,A)$ is a relative CW-complex which is $(l+1)$-connected. Then there exists a relative CW-complex $(Z,A')$ such that $A$ is a strong deformation retract of $A'$, $\cn A$ is a strong deformation retract of $Z$ and $(Z_{A'})^{l+1}=A'$. Let $r:A'\rightarrow A$ be the retraction and let $i_X:B\rightarrow X$ be the inclusion. Consider the pushout

\begin{displaymath}
\xymatrix{A' \ar[r]^(.5){gr} \ar[d]_{i_{A'}} \ar@{}[rd]|{push} & B \ar[d]^{i_Y} \\ Z \ar[r]_(.5){f'} & Y}
\end{displaymath}

Then $(Y,B)$ is a relative CW-complex with $(Y_B)^{l+1}=B$, and hence it is $(l+1)$-connected. The inclusions $i:A\rightarrow A'$ and $j:\cn A\rightarrow Z$ and the identity map of $B$ induce a map $\varphi:X\rightarrow Y$ with $\varphi i_X=i_Y \id_B$. Now, $i_A$, $i_{A'}$ are closed cofibrations and $i$, $j$ and $\id_B$ are homotopy equivalences, then, by proposition 7.5.7 of \cite{Br}, $\varphi$ is a homotopy equivalence. Thus, $(X,B)$ is $(l+1)$-connected.
\end{proof}

Note that the previous lemma can be applied when attaching a cell of any positive dimension, since attaching an $A$-n-cell is the same as attaching a $(\Sigma^{n-1}A)$-1-cell. The following lemma deals with the case in which we attach an $A$-0-cell. The proof is similar to the previous one.

\begin{lemma}
Let $A$ be an $l$-connected CW-complex, $B$ a topological space, and suppose $X$ is obtained from $B$ by attaching a 0-cell of type $A$ (i.e., $X=B\lor A$). Then $(X,B)$ is $l$-connected.
\end{lemma}

Now, using both lemmas we are able to prove the following proposition.

\begin{prop} \label{conndegofcwa}
Let $A$ be an $l$-connected CW-complex, and let $X$ be a CW($A$). Then the pair $(X,X^n)$ is $(n+l+1)$-connected.
\end{prop}

\begin{proof}
Let $r\leq n+l+1$ and $f:(D^r,S^{r-1}) \rightarrow (X^{n+1},X^n)$. We want to construct a map $f':(D^r,S^{r-1}) \rightarrow (X^{n+1},X^n)$ such that $f'(D^r)\subseteq X^n$, and $f\simeq f' \textrm{ rel }S^{r-1}$. Since $f(D^r)$ is compact, it intersects only a finite number of interiors of $(n+1)$-cells (note that $A$ is T1). By an inductive argument, we may suppose that we are attaching just one $(n+1)$-cell of type $A$, which is equivalent to attaching a $1$-cell of type $\Sigma^{n}A$. Since $\Sigma^{n}A$ is $(n+l)$-connected, $(X^{n+1},X^n)$ is $(n+l+1)$-connected. The result of the proposition follows.
\end{proof}

\begin{prop} \label{aconndegofcwa}
Let $A$ be an $l$-connected CW-complex, with $\dim(A)=k \in\mathbb{N}_0$, and let $X$ be a CW($A$). Then the pair $(X,X^n)$ is $A$-$(n-k+l+1)$-connected.
\end{prop}

\begin{proof}
We prove first the $A$-$0$-connectedness in case $k\leq n+l+1$. We have to find a dotted arrow in a diagram
\begin{displaymath}
\xymatrix{\ast \ar[r] \ar[d] \ar@{}[rd]|(.65){\simeq}|(.35){\mid\parallel} & X^n \ar[d]^i \\ A \ar[r]_f \ar@{.>}[ru] & X}
\end{displaymath}
This map exists because $A$ is a CW-complex with $\dim(A)=k$ and $(X,X^n)$ is $(n+l+1)$-connected.

Now we prove the $A$-$r$-connectedness in case $1 \leq r\leq n-k+l+1$. By lemma \ref{relhomotlemma}, it suffices to find a dotted arrow in a diagram
\begin{displaymath}
\xymatrix{\Sigma^{r-1}A \ar[r] \ar[d] \ar@{}[rd]|(.65){\simeq}|(.35){\mid\parallel} & X^n \ar[d]^i \\ \cn\Sigma^{r-1}A \ar[r]_(.6){f} \ar@{.>}[ru] & X}
\end{displaymath}
This map exists because $(\cn\Sigma^{r-1}A,\Sigma^{r-1}A)$ is a CW-complex of dimension $r+k$, $(X,X^n)$ is $(n+l+1)$-connected, and $r+k\leq n+l+1$.
\end{proof}

\email{ gminian$@$dm.uba.ar, emottina$@$dm.uba.ar}
\end{document}